\documentstyle[leqno,amssymb]{article}

\def\Z{\Bbb Z}
\def\N{\Bbb N}
\def\R{\Bbb R}
\def\CC{{\mathcal C}}
\def\IC{\Bbb C}

\def\endproof{\unskip\nobreak\kern5pt\nobreak\vrule height4pt width4pt depth0pt
\vskip4pt plus2pt}
\def\sq{\unskip\nobreak\kern5pt\nobreak\vrule height4pt width4pt depth0pt}

\def\endproof{\vrule height 0.5em depth 0.2em width 0.5em}
\newbox\tbox
\newbox\aubox
\newbox\adbox
\def\title#1{\setbox\tbox=\hbox{\let\\=\cr
\baselineskip14pt\vbox{\Large\bf\tabskip 0pt plus15cc
\halign to\hsize{\hfil\ignorespaces \uppercase{##}\hfil\cr#1\cr}}}}
\newbox\abbox
\setbox\abbox=\vbox{\vglue18pt}
\def\author#1{\setbox\aubox=\hbox{\let\\=\cr
\baselineskip12pt\vbox{\tabskip 0pt plus15cc
\halign to\hsize{\hfil\ignorespaces \uppercase{{##}}\hfil\cr#1\cr}}}%
\global\setbox\abbox=\vbox{\unvbox\abbox\box\aubox\vskip8pt}}
\def\address#1{\setbox\adbox=\hbox{\let\\=\cr
\baselineskip12pt\vbox{\it\tabskip 0pt plus15cc
\halign to\hsize{\hfil\ignorespaces {##}\hfil\cr#1\cr}}}%
\global\setbox\abbox=\vbox{\unvbox\abbox\box\adbox\vskip16pt}}
\def\makemytitle{
\begin{center}
\box\tbox
\box\abbox
\end{center}}

\begin{document}
%%%%%%%%%%%%%%%%%%%%%%%%%%%%%%%%%%%%%%%%%
\title{The Garden of Quantum Spheres}
\author{Ludwik D\c{a}browski}
\address{Scuola Internazionale Superiore di Studi Avanzati,\\
Via Beirut 2-4, I-34014, Trieste, Italy}
\makemytitle
%\vspace{0.5cm}
\abstract{
\begin{center}
A list of known quantum spheres 
of dimension one, two and three is presented.
\end{center}
}
%%%%%%%%%%%%%%%%%%%%%%%%%%%%%%%%
%%%%% remove for bcp:
\vspace{0.5cm}
\vspace{4mm}
\begin{center}
{\small
\noindent {\it M.S.C.}: 81R60, 81R50, 20G42, 58B34, 58B32, 17B37. \\
\noindent {\it Key words and phrases}: Noncommutative geometry, quantum spheres.}\\
\vspace{4mm}
Ref. SISSA 79/2001/FM\\
\end{center}
%%%%%%%%%%%%%%%%%%%%%%%%%%%%%%%%%%%%%%%%%
%\keywords{Noncommutative geometry, quantum spheres}
%\mathclass{Primary 58B34; Secondary 17B37.}
%\thanks{The author thanks the organizers for the invitation
%and `Geometric Analysis' Research Training Network HPRN-CT-1999-00118 of E. C.
%for the financial support.}
%\abbrevauthors{L.\ D\c abrowski}
%\abbrevtitle{Quantum Spheres}

\vspace{0.6cm}
\section*{0. Introduction.}~\\
\indent
Recently, examples of quantum spheres cropped in abundance in the literature.
The goal of this note is to aid the book-keeping of these newly emerged species 
by systematically comparing their basic properties.

As is customary in noncommutative geometry, these quantum spaces are described
and studied in terms of certain noncommutative algebras, 
generalizing the usual correspondence between spaces and function algebras. 
Here I am concerned mainly with `deformations' of the $*$-algebra 
of polynomials on the sphere $S^n$ and their enveloping $C^*$-algebras.
The $*$-algebras are usually given in terms of generators and relations.
Some of these relations can be regarded as deformations 
of the commutation relations and some as 
deformations of the sphere relation $\sum_j^{n+1} x_j^2 = 1$.
The classical spheres are often particular members of the family, or `limit' cases.
Here, by the dimension of such a quantum sphere I understand just the number $n$.

In this note the examples in lowest dimensions (one, two and three)
are listed (in alphabetic order), 
which appeared in the literature known to the author, 
without pretending to be complete or exhaustive.
Most of them have a $C^*$-algebraic version and often the deformation
forms a continuous field of $C^*$-algebras. 
The smooth structure is described only in a few cases.
Trying to uniformize the notation somewhat I mention 
some of their properties such as classical subspaces,
($C^*$-algebraic) K-groups and other remarks. 

It turns out that among basic building blocks of quantum spheres are
noncommutative tori and discs.
The $C^*$-algebra of the noncommutative torus $T_\theta$, $0\leq \theta <1$, 
is generated by two unitaries $U, V$ with the relation 
$UV = e^{2\pi i\theta} VU$~ \cite{r1}.
The $C^*$-algebra of a quantum disc $D_{\mu, q}$,
$0\leq \mu\leq 1$, $0<q\leq 1$, is generated by $z$ with the relation 
$qzz^* -z^*z = q-1 +\mu (1-zz^*)(1-z^*z)$~ \cite{kl}.
For $0\leq\mu<1-q$ and for $q=1$, $0<\mu <1$, they are 
known \cite{kl,s1} (see also \cite{hms1} and references therein) 
to be all isomorphic to the Toeplitz algebra \mbox{${\cal T}$} and, in turn, 
also to the $C^*$-algebra generated by the raising (shift) operator
$e_j\mapsto e_{j+1}$ on $\ell^2({\N})$, \cite{s1}.
% p. 222. 

%Throughout this note I denote natural numbers as $\N$, integers as $\Z$,
%real numbers as $\R$, 
In the sequel I denote the algebra unit as $I$
and the compact operators on a 
(infinite dimensional, separable) Hilbert space as $\cal K$.
%\newpage
\section*{1. Quantum circle.}
\subsection*{1.1. Markov \cite{m2}.}~\\
\underline{Generators}:~ $A, B$.~~
\underline{Relations}: 
$$
A=A^*, ~~  B=B^*, ~~A^2 + B^2 = I.
%, ~~AB = q BA
\nonumber 
$$
\underline{Classical subset}: $S^1$. 
~~\underline{K-groups}: $K_0={\Z}$, $K_1=0$ \cite{n1}.
~~\underline{Remarks}:
\begin{enumerate}
\item 
The universal (unital) $C^*$-algebra for these relations 
can be represented as the algebra of all
continuous $2\!\times\!2$  ~matrix-valued functions on the rectangle 
$[0, \pi/2 ]\times [0, \pi]$ satisfying certain boundary conditions.
\item  
The isomorphic universal unital $C^*$-algebra of the 
quantum disc at $q=-1, \mu = 0$ 
(generated by $z$ with the relation $zz^* + z^* z = 2I$)
has been independently studied in \cite{nn}.
\item 
Note that this example does not really fit to our list (of deformations).
% of quantum spheres $*$-algebras).
However at least it could be 
%To fit better our classification scheme the relations above could be (?)
supplemented by a commutation relation, e.g. $AB = q BA$, with $q= \pm 1$.
Then the case $q=1$ corresponds to the classical $S^1$ while 
in the case $q= -1$ the $C^*$-algebra is isomorphic to
the universal $C^*$-algebra of the free product of groups $Z_2 * Z_2$
with $K_0={\Z}^3$, $K_1=0$ and the 
classical subspace consisting of four points.
\end{enumerate}

\section*{2. Quantum 2-spheres.}
\subsection*{2.1. Bratteli,~Elliott,~Evans,~Kishimoto \cite{beek1,beek2,bk1}.}~\\
\underline{Generators}:~ $A, B$. 
~~\underline{Parameter}:~ $-1< \lambda = \cos(2\pi\theta) < 1$, ~($0<\theta< 1/2 $).
~~\underline{Relations}: 
\begin{eqnarray}
A=A^* &,& B=B^* ,\nonumber \\
BABA &=& \left( 4\lambda^2 - 1 \right) ABAB - 
2 \lambda A^2 B^2 + 8\lambda \left( 1- \lambda^2 \right) (A^2 + B^2 - I),\nonumber\\
A^2B + BA^2 &=& 2 \lambda ABA + 4 (1- \lambda^2 ) B,\nonumber \\
AB^2 + B^2 A &=& 2\lambda BAB + 4 (1- \lambda^2 ) A. \nonumber
\end{eqnarray}
\underline{Classical subset}: $\emptyset$. ~~\underline{K-groups}: 
$K_0={\Z}^6$, $K_1=0$  \cite{n3,k}.
~~\underline{Remarks}:
\begin{enumerate}
\item 
Has been introduced, via $A= U+ U^{-1}, V+V^{-1}$,
as a fixed point algebra of the `flip' automorphism 
$\sigma: U\mapsto U^{-1}, V\mapsto V^{-1}$
of the noncommutative torus $T_\theta$.
\item  
Classically (for $\theta =0$) this is a `pillow'
(a smooth 2-sphere with four corners).
After the deformation (for $0<\theta < 1/2$) this geometry 
manifests in $K_0={\Z}^6$ (with four generators besides 
$I$ and the Bott projector).
\item  
For $\lambda =1$ ($\theta =0$) this $C^*$-algebra 
indeed corresponds to the classical $S^2$ as it turns out that
in this case either of the last two relations implies that $AB=BA$.
\item  
For irrational $\theta $ these $C^*$-algebras are simple, 
approximately finite-dimensional and with a unique trace state.
%(is a corner in) Bratteli.~et.~al./Kumjian/Walters Example 2.2.
\item 
A closely related (strongly Morita equivalent for $\theta\neq 1/2$)
is the crossproduct $C^*$-algebra 
$T_\theta {>\kern-5pt\triangleleft}_\sigma {\Z}_2$,
where the generator of ${\Z}_2$ acts
on the noncommutative torus $T_\theta$ by the `flip' automorphism 
$\sigma: U\mapsto U^{-1}, V\mapsto V^{-1}$, cf. \cite{beek1,beek2,bk1,k,w1}.
It is generated by three unitaries $U, V, W$ with the relations
$VU = e^{2\pi i{\theta}} UV$, $WUW = U^*$, $WVW = V^*$, $W^2 = I$.
For irrational $\theta, \theta^\prime $ these 
$C^*$-algebras are isomorphic iff $\theta = \theta^\prime $ or
$\theta = 1-\theta^\prime $.
For rational $\theta = p/q, \theta^\prime =p^\prime / q^\prime$,
with $p, p^\prime$ and also $q, q^\prime$ relatively prime,
they are isomorphic iff $q = q^\prime $.
\item 
It admits other presentations \cite{k}
with three selfadjoint unitaries $X, Y, Z$
and the relation $XYZ = e^{2\pi i{\theta}} ZYX$,
or with four selfadjoint unitaries $X, Y, Z, T$
and the relation $XY = e^{\pi i{\theta}} TZ$.
\end{enumerate}

\subsection*{2.2. Calow,~Matthes \cite{cm1}.}~\\
\underline{Generators}:~ $A, B$.
~~\underline{Parameters}:~ $0<p, q<1$.
~~\underline{Relations}: 
\begin{eqnarray}
A &=&A^*, \nonumber \\
B^*B-qBB^* &=& (p\! -\! q)A+ I-p, \nonumber \\
AB-pBA &=& (1\! -\! p)B \nonumber \\
(I-A)(BB^*-A) &=& 0.\nonumber 
\end{eqnarray}
\underline{Classical subset}: ~$S^1$.
~~\underline{K-groups}: ~$K_0={\Z}^2$, $K_1= 0$ \cite{cm1}.
~~\underline{Remarks}:
\begin{enumerate}
\item 
Obtained (at the $*$-algebra level) by glueing the quantum disc 
$D_{0, p}$ with $D_{0, q}$.
\item 
Classically (for $p=1=q$) this $*$-algebra describes a closed cone
with one vortex and one circular edge.
\item 
It is non $*$-isomorphic \cite{cm1}
with any of the Podle\'s quantum spheres.
However, its universal $C^*$-algebra is isomorphic 
to the $C^*$-algebra of the generic ($s>0$) Podle\'s quantum sphere,
and as such can be realized as the Cuntz-Krieger algebra of a certain graph 
or as a quantum double suspension of two points \cite{hs}.
\end{enumerate}

\subsection*{2.3. Gurevich,~Leclercq,~Saponov \cite{gls}.}~\\ 
\underline{Generators}:~ $A, B$.
~~\underline{Parameters}:~ $h \geq 0$, $ q > 0$.
~~\underline{Relations}: 
$$
A = A^*, ~~qB^*B + q^{-1}BB^* +q^2(q+ q^{-1}) A^2 = q+ q^{-1},\nonumber
$$
$$
q^2 AB -BA = hB, ~~BB^* - B^*B = (1- q^{4}) A^2 + h(1+q^2)A.\nonumber
$$
\underline{Classical subset}: $S^2$ for $h=0, q=1$; $S^1$ for $h=0, q\neq 1$; 
a point for $\pm (q^4-1)/q = h \neq 0$ and $\emptyset$ otherwise.
~~\underline{Remarks}:
\begin{enumerate}
\item 
I have set the (radius) parameter $-\alpha =1$ in \cite{gls} 
by overall rescaling of the generators,
which relate to those used in \cite{gls} as $B=b$ and $A=g/(1+q^2)$,
and the parameters as $h= q^{-1}\hbar$ ($q$ is unchanged). 
\item 
This example deforms the universal enveloping algebra of $su(2)$
with a constrained value of the quadratic Casimir element.
\item 
The $C^*$-algebra completion is in general not possible
but can be accomplished in particular cases.
\item 
The one-dimensional subfamily with $h=0$ coincides
with the one parameter subfamily of equatorial quantum 
Podle\'s \cite{p} spheres.
\item 
Another particular one-dimensional subfamily with $q=1$ coincides, 
using the variables $Z= B/\sqrt{h^2+4}$, $H= A/\sqrt{h^2+4}$ and 
$\mu= h/\sqrt{h^2+4}$ (when treated as a formal parameter)
with the quantum 2-sphere of Omori,~Maeda,~Miyazaki,~Yoshioka \cite{ommy}.
Any member of this subfamily 
forms a basis of a quantum principal $U(1)$-bundle 
with a total space being a formal deformation of $S^3$
(with a noncentral formal parameter $\mu$), c.f. Section 4.
\end{enumerate}

\subsection*{2.4. Natsume \cite{n2}.}~\\
% or Natsume, Olsen \cite{no2}.}~\\ ?????????????
\underline{Generators}:~ $A, B$.
~~\underline{Parameter}:~ $t \in {\R}$.
~~\underline{Relations}: 
$$
A=A^*, ~
B^* B + A^2 = I, 
~BB^* + (t BB^* + A)^2 =I,
~BA - AB = t BB^*B.$$
\underline{Classical subset}: two points for $t\neq 0$.
~~\underline{K-groups}: $K_0={\Z}^2$, $K_1= 0$ \cite{no2}.\\
\underline{Remarks}:
\begin{enumerate}
\item 
Motivated by Poisson geometry.
\item 
For $t\in [0, 1/2[$ the enveloping $C^*$-algebras are of type I 
and are isomorphic \cite{no2}
to certain extension of ${\Bbb C}^2$ by the crossproduct $C^*$-algebra 
$C_0\left(\, ]-\! 1,1[\, \right){>\kern-5pt\triangleleft}_{\alpha_t}{\Z}$,
or equivalently by ${\cal K}\otimes C(S^1)$
(the generator of ${\Z}$ acts by automorphism $\alpha_t$,
given by a pull back of the homeomorphism 
$x\mapsto tx^2+x-t$, $\forall x\in\, ]-\! 1, 1[$~, 
which is topologically conjugate to the translation by 1 on ${\Bbb R}$)
\cite{no2}.
\item 
They form a continuous field of $C^*$-algebras over $[0, 1/2[$,
which is trivial over $]0, 1/2[$ (in particular they are all isomorphic
for $t\in\, ]0, 1/2[$~).
They also constitute a strong deformation of $C(S^2)$ \cite{no2}.
\end{enumerate}

\subsection*{2.5. Podle\'s \cite{p}.}~\\
\underline{Generators}:~ $A, B$.~~
\underline{Parameters}:~ $0\leq q <1$, ~$ 0\leq s \leq 1$. ~~
\underline{Relations}:
$$
A^* = A, ~~BA=q^2AB, \nonumber
$$
$$
BB^*=-q^4A^2+(1-s^2)q^2 A +s^2 I,
~~B^*B= -A^2 + (1-s^2)A + s^2 I. \nonumber
$$
\underline{Classical subset}: 
a point if $c\! =\! 0$; $S^1$ if $c\!\in\, ]0, \infty]$.
~\underline{K-groups}: $K_0={\Z}^2$,
$K_1=0$ \cite{mnw2}.\\
\underline{Remarks}:
\begin{enumerate}
%{\hspace*{-0.5cm}
\item
%i) Also $q=0$ can be allowed.
Discovered as homogeneous $SU_q(2)$-spaces.
\item 
In order to write the relations in a uniform way 
I parametrize the whole family as in \cite{hms4}
by the parameter $0\leq s \leq 1$, 
related to the parameter $ 0\leq c \leq \infty $ in \cite{p} by 
$s=2\sqrt{c}/(1+\sqrt{1+4c})$ (and $c=(s^{-1}-s)^{-2}$). 
Also, I use the Podle\'s generators $A, B$ rescaled, 
iff~ $0\leq s < 1$ (i.e., $ 0\leq c < \infty $), 
by $1-s^2$.
\item 
%\hspace*{-0.5cm}
Any memebers of this family describes a `round' quantum sphere, 
in the sense that the Cartesian coordinates can be found 
i.e., three selfadjoint elements which generate the $*$-algebra, 
commute among themselves for $q=1, s=0$ and whose squares sum up to $I$.
\item 
The case $s=0$ (i.e., $c=0$) when the last two relations read
$B B^* = q^2 A -q^4 A^2$ and $B^* B = A - A^2$,
is known as the {\em standard Podle\'s sphere}.
It can be viewed as a quotient sphere $SU_{q}(2)/{U(1)}$ 
in the spirit of the Hopf fibration (c.f. Example 3.4).
For all $0\leq q <1$, the corresponding $C^*$-algebras are isomorphic 
to the minimal unitization of the compacts $\cal K$. 
\item 
For $0< s \leq 1$ (i.e., $0 < c \leq \infty $) 
the related $C^*$-algebras are all isomorphic \cite{s1}
to certain extension of $C(S^1)$ by $ \cal K\oplus \cal K$,
or extension of $\cal T$ by $ \cal K$.
They can also be described, at the $C^*$-algebra level, 
as two quantum discs glued along their boundaries $S^1$ \cite{s1}, 
see also \cite{bk2,cm1};
%$D_{\mu, 1}$ [8], or by glueing the quantum disc $D_{0, p}$ with $D_{0, q}$ \cite{cm1};
as the Cuntz-Krieger algebra of a certain graph 
or as a quantum double suspension of two points \cite{hs}.
\item 
The case $s=1$ (i.e. $c=\infty$) when the last two relations read
$ B B^* = -q^4 A^2 + I$ and $B^* B = - A^2 + I $,
is known as the {\em equatorial Podle\'s sphere}.
It is easily seen to be $*$-isomorphic to the two dimensional
Euclidean sphere, introduced in \cite{frt}.
As such, it admits a higher (even) dimensional generalization.
Also, it contains the $*$-algebra of quantum disk,
which can be geometrically interpreted as collapsing this quantum 2-sphere
by the reflection with respect to the equatorial plane \cite{hms1}.
Moreover, it is isomorphic to the quotient of the underlying 
$*$-algebra of Example 3.4. (with parameter $q^2$) by the relation $b=b^*$.
The geometric meaning of this is that the equatorial Podle\'s sphere
embeds as an equator in $SU_{q^2}(2)$ thought of as a quantum 3-sphere \cite{hms3}.
Hence for fixed $q$,
the path $0\leq s\leq 1$ of Podle\'s spheres 
can be viewed as an interpolation between the quotient sphere
$SU_{q}(2)/{U(1)}$ and the embedded (equator) 2-sphere in $SU_{q^2}(2)$. 
\end{enumerate}

\section*{3. Quantum 3-spheres.}
\subsection*{3.1. Calow,~Matthes \cite{cm2}.}~\\
\underline{Generators}:~ $a, b$.~~
\underline{Parameters}:~ $0<p, q<1$.~~
\underline{Relations}: 
$$
{a}^*{a} - q{a} {a}^* = 1-q,\nonumber \\
~{b}^* {b}-p{b} {b}^* = 1-p, \nonumber\\
~ab=ba,~a^*b=ba^*, \nonumber\\
(I-{a}{a}^*)(I-{b} {b}^*)= 0.\nonumber 
$$
\underline{Classical subset}: ~$S^1\times S^1$.
~~\underline{K-groups}: ~$K_0={\Z}$, ~$K_1={\Z}$ \cite{hms2}.
~~\underline{Remarks}:
\begin{enumerate}
\item 
As a $*$-algebra obtained by glueing the quantum solid torus $D_{0, p}\times S^1$ with
$D_{0, q}\times S^1$.
\item 
As a $C^*$-algebra isomorphic to 
$ (\cal T\otimes \cal T)/(\cal K\otimes \cal K)$,
or also to the Cuntz-Krieger algebra of a certain graph \cite{hms2}.
\item 
Forms a locally trivial, globally nontrivial \cite{cm2},
in fact noncleft \cite{hms2},
quantum principal $U(1)$-bundle (Hopf-Galois extension) over 
the quantum $S^2$
%2-sphere 
of Example 2.2.
\end{enumerate}

\subsection*{3.2. Connes,~Dubois-Violette \cite{cdv}.}~\\
\noindent 
\underline{Generators}:~ $x^0, x^1, x^2, x^3$. 
~~\underline{Parameters}:~ $\pi > \varphi_{1}\geq\varphi_{2}\geq\varphi_{3}\geq 0$.
~~\underline{Relations}: 
\begin{eqnarray}
x^0 = (x^0)^*,~~~~~~~~ x^1 = (x^1)^*,&& x^2= (x^2)^*,~~~~~~~~ x^3= (x^3)^*,\nonumber \\
\cos(\varphi_{1})(x^0x^1-x^1x^0) & = 
&i\sin(\varphi_{2}-\varphi_{3})(x^2x^3+x^3x^2),\nonumber\\
\cos(\varphi_{2})(x^0x^2-x^2x^0)& = &i\sin(\varphi_{3}-\varphi_{1})(x^3x^1+x^1x^3),\nonumber\\
\cos(\varphi_{3})(x^0x^3-x^3x^0)& = &i\sin(\varphi_{1}-\varphi_{2})(x^1x^2+x^2x^1),\nonumber\\
\cos(\varphi_{2}-\varphi_{3})(x^2x^3-x^3x^2)&= &-i\sin(\varphi_{1})(x^0x^1+x^1x^0),\nonumber\\
\cos(\varphi_{3}-\varphi_{1})(x^3x^1-x^1x^3)&= &-i\sin(\varphi_{2})(x^0x^2+x^2x^0),\nonumber\\
\cos(\varphi_{1}-\varphi_{2})(x^1x^2-x^2x^1)&= &-i\sin(\varphi_{3})(x^0x^3+x^3x^0),\nonumber\\
(x^0)^2+(x^1)^2+(x^2)^2+(x^3)^2 &=& 1. \nonumber
\end{eqnarray}
\underline{Classical subset}: generically discrete.
~~\underline{K-groups}: will be studied in part II of \cite{cdv}.
~~\underline{Remarks}:
\begin{enumerate}
\item 
Obtained as a unital  $*$-algebra generated by four generators $u_{jk}$,
$j,k\in \{1,2\}$ such that $u$ as a $2\times 2$ matrix is unitary
and $ch_{1/2} (u):= \sum_{j,k} (u_{jk}u^*_{kj}-u^*_{jk}u_{kj})=0$.
\item 
The particular one parameter subfamily 
$\varphi_{1}=\varphi_{2}=-\frac{1}{2}\theta$ and 
$\varphi_{3}=0$ coincides, 
using the variables $Z=x^0 +i x^3, W=x^1 +i x^2$, with the particular 
one parameter subfamily of Matsumoto quantum spheres 3.3,
when $\Theta = \theta$ is a constant function,
and thus also with the Natsume, Olsen \cite{no} family.
It fulfills all the properties of 
a noncommutative manifold in the sense of \cite{c1}
and has a higher (odd) dimensional generalization.
\end{enumerate}

\subsection*{3.3. Matsumoto \cite{m3,m4}.}~\\
\underline{Generators}:~ a pair $Z, W$ of normal operators on a Hilbert space.\\
\underline{Parameters}:~ real valued continuous functions $\Theta$ on the closed
interval $[0, 1]$.\\
~~\underline{Relations}: 
$$
ZW = e^{2\pi i{\widehat\Theta}(Z^*Z)} WZ\ , ~~Z^* Z + W^* W =I \ ,
$$
where ${\widehat\Theta}(Z^*Z)$ 
stands for the self-adjoint operator obtained by the functional calculus 
of the operator $Z^*Z$ from the continuous function $\Theta$.\\
\underline{Classical subset}: 
$\left((\Theta^{-1}({\Z})\cap \; ]0,1[\; )\times T^2\right)
\cup\left((\Theta^{-1}({\Z})\cap\{0,1\})\times S^1\right)$.\\
~~\underline{K-groups}: ~$K_0={\Z}$, ~$K_1={\Z}$.
~~\underline{Remarks}:
\begin{enumerate}
\item 
Obtained by glueing two quantum solid tori described by 
crossproduct $C^*$-algebra 
$\CC (D) {>\kern-5pt\triangleleft}_\Theta \Z$,
where the generator of $Z$ acts on the 2-disk $D$ as a rotation by angle 
$\Theta (r)$, and $r\in [0, 1]$ is the radial coordinate on $D$.
%, constructed using $T_\theta$ as a building block.
\item 
Forms a quantum principal $U(1)$-Hopf fibration
over the usual 2-sphere.
\item 
When the function $\Theta$ is a constant number $\theta$,
the $C^*$-algebra generated by the relations above 
has been studied by Natsume and Olsen \cite{no}
and shown to be isomorphic to the universal $C^*$-algebra generated 
by two normal operators $T, S$ satisfying 
$TS = e^{2\pi i{\theta}} ST$, $(I-T^*T)(I-S^*S)=0$ and $\| T\| = 1= \| S\|$, 
introduced in \cite{m3}.
Its classical subset is $S^1\sqcup S^1$
%its K-groups are $K_0={\Z}$, $K_1={\Z}$, 
and it has 
%$2n\! -\!1 $
odd-dimensional generalization \cite{no}.
It coincides with the particular 
one parameter subfamily 
$\varphi_{1}=\varphi_{2}=-\frac{1}{2}\theta$ and 
$\varphi_{3}=0$ of Connes,~Dubois-Violette \cite{cdv}.
\end{enumerate}

\subsection*{3.4. Woronowicz \cite{w2}.}~\\
\underline{Generators}:~ $a, b$.
~~\underline{Parameters}:~ $ q\in \IC$.  
~~\underline{Relations}: 
$$ba = q ab, ~a^*b = q ba^*,
%b^*a = qab^*, a^*b^*=q b^* a^*,
~aa^* + bb^*=I, ~a^*a + |q|^2 bb^* =I, ~bb^*=b^*b.$$
\underline{Classical subset}: ~$S^1$ for $|q|\neq 1$, 
$S^1\sqcup S^1$ for $|q|=1$ and $q\neq 1$, $S^3$ for $q=1$.\\
~~\underline{K-groups}: ~$K_0={\Z}$, ~$K_1={\Z}$ for $q>0$ \cite{mnw1}.
~~\underline{Remarks}:
\begin{enumerate}
\item 
Has been discovered for $-1 \leq q \leq 1$, $q\neq 0$, 
as a family of quantum groups $SU_q(2)$. 
\item 
Here I generalize the range of the parameter to $ q\in \IC$,
this `interpolates' between the original Woronowicz family $q\in \R$,
and the one-parameter subfamily $|q|=1$ which coincides (if $q=e^{i\theta}$) 
with Natsume and Olsen \cite{no} family and also with the particular 
one parameter subfamily $\varphi_{1}=\varphi_{2}=-\frac{1}{2}\theta$ and 
$\varphi_{3}=0$ of the example 3.2 of Connes,~Dubois-Violette \cite{cdv}.
%>>>
The transformation $a\mapsto a^*$,  $b\mapsto -q b$, and $q\mapsto 1/q$ 
defines a $*$-isomorphism for $q\neq 0$, hence it suffices to restrict
to the range $|q|\leq 1$.
%The one-parameter subfamily $0<q<1$ is 
%the original Woronowicz $q\in \R$ 
%$*$-algebra when 
\item 
For $0 \leq q < 1$ the members of this family are easily seen to be 
$*$-isomorphic to the three dimensional Euclidean spheres
introduced in \cite{frt} and also to the three dimensional unitary spheres
introduced as quantum homogeneous spaces of $SU_q(n)$ in \cite{vs}.
It turns out that also their higher dimensional generalizations,
the quantum Euclidean spheres and the quantum unitary spheres,
are $*$-isomorphic at any given odd dimension.\footnote{
To our knowledge this simple fact, which was observed during a conversation 
with G.\ Landi, E.\ Hawkins and F.\ Bonechi, has not been presented before.}
\item 
For all $ q\in \IC$  these $*$-algebras have a $C^*$-algebraic
version. For $q=1$ this is just $C(S^3)$. The case $q=-1$ has been studied
in \cite{zak}.
For $0 \leq q < 1$ these $C^*$-algebras are all isomorphic
to certain extension of $C(S^1)$ by $C(S^1)\otimes \cal K$,
which can also be described as the Cuntz-Krieger algebra of certain graph 
\cite{hs} or as a quantum double suspension of the circle \cite{hs}.
\item
 For $0 < q \leq 1$ forms a quantum principal $U(1)$-fibre bundle 
over the Podle\'s quantum 2-sphere of Example 2.5 in the sense of 
Hopf-Galois extensions (if $s=0$) or coalgebra-Galois extensions
(if $s\in\, ]0,1]$)
\cite{bm,hm,b}.
\item 
For some further noncommutative-geometric aspects 
see \cite{c2} and references therein.
\end{enumerate}

\section*{4. Final comments.}
~\\
Some finite dimensional algebras (in a sense corresponding
to zero dimensional quantum spaces, which nevertheless possess
certain properties of 2-spheres) have also been studied.   
For instance the classification \cite{p} of $SU_q(2)$-homogeneous spaces,
besides the family 2.5 of Podle\'s quantum spheres, 
also includes a discrete series of full 
%$N\times N$ 
matrix algebras ${\rm Mat}_N$.
It has been observed in \cite{gms} that this family of `quantum spheres' can be 
equipped with an additional structure, notably a sequence of injections 
${\rm Mat}_N \to {\rm Mat}_{N+1}$,
which are morphisms in the category of $U_q(su(2))$-modules.
For $q=1$, this agrees with the {\em fuzzy-sphere} philosophy of \cite{m1}.
Therein, the $N\times N$ matrix algebras are considered
as $U(su(2))$-modules together with $U(su(2))$-module injections
that form a direct system whose limit is the algebra of polynomials on $S^2$ \cite{h}.
One can also show that the matrix algebras converge 
to the sphere for quantum Gromov--Hausdorff distance \cite{r2}.
%math.OA/0108005
Furthermore, these matrix algebras can be viewed as representations 
of the universal enveloping algebra of $sl(2)$ 
with the value $\frac{N^2-1}{4}$ of the quadratic Casimir element.
This is why the discrete family of Podle\'s spheres
can be thought of as the family of q-fuzzy spheres \cite{gms}
(the aforementioned injections are $U_q(su(2))$-linear).

There are other examples of quantum spheres which do not fit 
exactly to our lists as they are not deformations in our sense 
(families of $*$-algebras). 
In \cite{ommy,omym} a formal deformation of $S^3$ (as a contact manifold)
is provided with a invertible non-central deformation `parameter' $\mu$,
generators $a, b$ and relations 
$$
\mu = \mu^*, ~a^*a + b^*b= I, \mu^{-1} a - a\mu^{-1}  = -a , 
~\mu^{-1} b - b \mu^{-1}  = -b,
$$
$$
~ba = ab, ~ab^* = (1-\mu ) b^*a, 
aa^* - (1-\mu )a^* a=\mu,  ~bb^* - (1-\mu )b^*b=\mu .
$$
This deformation yields a certain `smooth' algebra admitting a $U(1)$-action
that is principal in the sense of Hopf-Galois theory \cite{bdz}
(the base space of the pricipal fibre bundle is given by the quantum 2-sphere of 
of Example 2.3.5.)
Another example of a noncommutative Hopf fibration given by 
the principal $U(1)$-action on a super 3-sphere was studied in \cite{dgh}.
There are also examples of quantum complex spheres related to
the Jordanian quantum group $SL_h(2)$ \cite{cs,z}.

As far as four-dimensional quantum spheres are concerned,
recently several examples together with instanton bundles over them,
have been constructed.
They indicate a wealth even greater than that of the known 
one, two and three dimensional examples.
It should be mentioned that in general various principles 
were employed to proliferate the examples of quantum spheres,
such as Poisson, contact or homogeneous structure. 
We have encountered also several types of glueing and quotienting.
One of the tools is the usual suspension operation 
which can be used to find links between different examples
and to produce new samples of one dimension greater.
Note for instance that the suspension of Example 3.2.2 
is just the Connes,~Landi noncommutative $4$-sphere 
\cite{cl},
while the suspension of Example 3.4 occurs in 
\cite{dl},
c.f. also 
\cite{dlm}.
A $C^*$-algebraic noncommutative double suspension has been also 
mentioned in examples 2.2, 2.5 and 3.4.
A kind of quantum double suspension at the $*$-algebra level, 
which raises the dimension by two and yields different families 
of quantum 4-spheres, has been employed in \cite{s2} and \cite{bg},
The example of a quantum 4-sphere presented in \cite{bct}, 
with its classical subspace being just a point,
is yet another type of double suspension \cite{bct2}
of the standard Podle\'s quantum sphere,
motivated by Poisson structure.
Another method to obtain more examples employs the {\em twisting} 
(see e.g., \cite{cl,cdv}).
%a different nature

However, it seems that it is quite premature yet to attempt any 
classification of quantum 4-spheres 
and that is certainly beyond the scope of this note.

\vspace{5mm}
{\small \noindent{\bf Acknowledgments:} 
The author thanks the organizers for the invitation,
the `Geometric Analysis' Research Training Network HPRN-CT-1999-00118 of E. C.
for the financial support,
and E.~Hawkins and  P.~M.~Hajac for several comments on the ma\-nu\-script.}
%\vspace{5mm}

%\references{Nov}

%}
%\end{minipage}
\end{document}